\documentclass[a4paper,11pt]{article}
\usepackage{fullpage}
\usepackage{times}
\usepackage{subfigure}
\usepackage{vmargin}
\usepackage{fancyheadings}
\usepackage[cp1251]{inputenc}
\usepackage[english]{babel}
\usepackage{graphicx}
\usepackage{color}
\usepackage{amsmath,amssymb,amsthm,enumerate,bbm,color,verbatim,setspace}

\newtheorem{lemma}{Lemma}[section]

\newtheorem{Def}{Definition}[section]
\date{}
\title{Gamma-weibull kernel estimation of the heavy tailed densities}
\author{ L. Markovich$^{1}$\\
       $^{1}$ Institute of Control Sciences,
         Russian Academy of Sciences,\\
Moscow, Russia. \\
$^1$Corresponding author e-mail: kimo1@mail.ru}
\date{}
\begin{document}
\maketitle
\pagenumbering{arabic}
\begin{abstract}\noindent
We consider the nonparametric estimation
 of the univariate heavy tailed probability density function (pdf) with a support on
 $[0,\infty)$ by independent data. To this end we construct the new kernel estimator as a combination of the asymmetric
gamma  and weibull kernels, ss. gamma-weibull kernel. The gamma kernel is nonnegative, changes the shape depending on the position
on the semi-axis and possess good boundary properties for a wide class of densities. Thus, we use it to estimate the pdf  near the zero boundary. The weibull kernel is based on the weibull distribution which can be heavy tailed and hence we use it to estimate the tail of the unknown  pdf. The theoretical asymptotic properties  of the proposed density estimator like bias and variance are derived.
We obtain the optimal bandwidth selection for the estimate as a minimum of the mean integrated squared error (MISE). Optimal rate of convergence of the MISE for the density is found.
\end{abstract}
\medskip

\section{Introduction}
Estimation of the nonnegatively supported probability density functions (pdf) appear in many fields of the applied statistics. Such pdfs are used to model
a wide range of applications in engineering, signal processing  \cite{Dobrovidov:12}, medical research, quality control, actuarial science \cite{Furman,Hurlimann} and climatology \cite{Aksoy} among others.
\par The wide use of these pdfs in practice leads to the need of their estimation by finite data samples. One of the most common nonparametric pdf estimation methods are kernel estimators. However most of the known asymmetric estimators are oriented on light-tailed distributions. For example, for the univariate nonnegative independent identically distributed (iid) random variables (r.v.s), the  estimators with gamma kernels were proposed in \cite{Chen:20}.  The gamma kernel estimator was developed for univariate dependent data in \cite{TaufikBouezmarnia:Rom2}. In \cite{TaufikBouezmarnia:Rom} the gamma kernel estimator of the multivariate pdf for the nonnegative iid r.v.s was introduced. In \cite{Mar_2015} the gamma kernel estimator of the multivariate pdf and its gradient for the nonnegative dependent r.v.s was obtained. Other asymmetrical kernel estimators like inverse Gaussian and reciprocal inverse Gaussian estimators were studied in \cite{Scailet}. The comparison of these asymmetric kernels with the gamma kernel is given in \cite{BouSca}.
\par All these estimators are called
nonparametric since they do not require any preliminary information regarding
the parametric model of the distributions but only very common features
of the distributions like a number of continuous derivatives, for example.
The focus of our paper is on the nonparametric estimation of heavy-tailed densities
which are defined on a positive part
of the real axis. It is obvious, that the known classical estimators
cannot be directly applied to heavy-tailed distributions. These are characterized by slower
decay to zero of heavy tails than that of an exponential rate, the lack of
some or all moments of the distribution, and sparse observations at the tail
domain of the distribution.
\par The known approaches of the heavy-tailed density estimation are the kernel estimators with the heavy tailed kernels, the estimators based on the preliminary transform of the initial random variable
(rv) to a new one and "piecing-together approach" which provides
a certain parametric model for the tail of the density and a non-parametric
model to approximate the "body" of the density.
\par In this paper, we introduce the new kernel constructed from the gamma and the weibull kernel estimators. The new Gamma-weibull kernel has two smoothing parameters (bandwidths) and the third parameter - the width of the boundary domain of the gamma part of the kernel.
\par The paper is organized as follows. In Sec. \ref{sec:1} we introduce the new gamma-weibull kernel estimator. In Sec. \ref{sec:2} we obtain the bias and the variance of the pdf estimate. Using these results we derive the  optimal bandwidths and the corresponding rate of the optimal MISE.
\section{The Gamma-Weibull Kernel}\label{sec:1}
\par The term heavy-tailed is used to the class of probability density function whose tails are not exponentially bounded, i.e. there tails are heavier then the exponential pdfs tail \cite{Asmussen}. Some authors define the heavy-tail pdfs as the pdf with some infinite power moments, for example variance. Let $\{X_i;i=1,2,\ldots\}$ be a strongly stationary sequence with an unknown pdf $f(x)$  and distribution function (df) $F(x)$ which are defined on the nonnegative semiaxes $x\in[0,\infty)$. In \cite{Embrechts} the following definition is given.
\begin{Def}The distribution of a r.v. $X$ is said to have the heavy-tail if
\begin{eqnarray*}\lim\limits_{x\rightarrow\infty}\textbf{P}\{X>x+y|X>x\}=\lim\limits_{x\rightarrow\infty}\overline{F}(x+y)/\overline{F}(x)=1,\quad y>0
\end{eqnarray*}
\end{Def}
The examples of such pdfs are Lognormal, Pareto, Burr, Cauchy, Weibull with with shape parameter less than 1 among others.
\par Our objective is to estimate the unknown pdf  by a known sequence of observations $\{X_i\}$. Since the pdf is asymmetric and
 can be heavy tailed we can not use the standard gaussian kernel estimator.
\par Let us construct the special kernel function which would be flexible on the domain near the zero boundary and could estimate the heavy tail of the distribution.
For the domain $x\in[0,a],a>0$ we use the non-symmetric gamma kernel estimator that was defined in \cite{Chen:20} by the formula
\begin{eqnarray*}
\widehat{f_G}_n(x) &=& \frac{1}{n}\sum\limits_{i=1}^{n}K_{\rho(x,h),\theta}(X_i)= \frac{1}{n}\sum\limits_{i=1}^{n}\frac{X_i^{\rho(x,h)-1}e^{-X_i/\theta}}{\theta^{\rho(x,h)}\Gamma(\rho(x,h))},\quad \rho,\theta>0.
\end{eqnarray*}
Here $\Gamma(\rho)$ is the gamma function evaluated at $\rho$ and $h$ is the bandwidth of the kernel. The shape parameters $\rho,\theta$ will be selected further.
\par For the domain $ x> a$ the Weibull kernel estimator is constructed
 \begin{eqnarray*}
 \widehat{f_W}_n(x) &=&\frac{1}{n}\sum\limits_{i=1}^{n}K_{k(x,b)}\left(X_i\right)= \frac{1}{n}\sum\limits_{i=1}^{n}\frac{k(x,b)}{\lambda}\left(\frac{X_i}{\lambda}\right)^{k(x,b)-1}\exp\left(-\left(\frac{X_i}{\lambda}\right)^{k(x,b)}\right),
 \end{eqnarray*}
 where the shape parameters are $\lambda>0$, $0<k<1$ and $b$ is the bandwidth of the kernel. Hence, the pdf estimator is the folowing
\begin{eqnarray}\label{12}
 \widehat{f_{GW}}_n(x)&=& \left\{
\begin{array}{ll}
\widehat{f_G}_n(x) &   \mbox{if}\qquad x\in
[0,a],
\\
\widehat{f_W}_n(x)& \mbox{if}\qquad x> a.
\end{array}
\right.
\end{eqnarray}
The latter kernel estimator has two bandwidth parameters $h$ and $b$ and one special parameter $a$.
The parameters $\rho(x,h)$, $k(x,b)$, $\lambda$ and $\theta$ can be found from the matching conditions
\begin{eqnarray}&&f_G(X,\rho(x,h),\theta)\Big|_{x=a}-f_W(a,k(x,b),\lambda)\Big|_{x=a}=0,\label{1}
\end{eqnarray}
\begin{eqnarray}
&&f'_G(X,\rho(x,h),\theta)\Big|_{x=a}-f'_W(a,k(x,b),\lambda)\Big|_{x=a}=0.\label{2}
\end{eqnarray}
From the condition \eqref{1} we can deduce that the shape parameters of the kernels are
\begin{eqnarray*}\rho(a,h)&=&k(a,b),\quad \theta=\lambda.
\end{eqnarray*}
From the condition \eqref{2} we can deduce that
\begin{eqnarray}\frac{\partial \rho(x,h)}{\partial x}\Big|_{x=a}&=&\frac{\partial k(x,b)}{\partial x}\Big|_{x=a}.
\end{eqnarray}
Hence, we can select any variety of $\rho(x,h)$ and $k(x,b)$ that satisfy the latter conditions to get some kernel estimators. Let us select for example the following parameters
\begin{eqnarray}\label{3}\rho(x,h)&=&\frac{x+c_1h}{a},\quad  k(x,b)=\frac{x+c_2b}{a}.
\end{eqnarray}
Hence, the bandwidth parameters satisfy the condition $h=bc_2/c_1$ when $x=a$. Since $k(x,b)<1$ for the heavy tailed Weibull pdf the parameters $c_1$, $c_2$ are some negative constants that we will select further.
\par As the measure of error of the proposed estimator \eqref{12} it is natural to consider the MISE which is defined as
\begin{eqnarray}\label{3a}
MISE(\widehat{f}(x)&=&\mathsf E\int\limits_0^\infty(f(x)-\widehat{f}(x))^2dx.
\end{eqnarray}
The unknown smoothing parameters $h$ and $b$ are obtained as the minima of \eqref{3a}.
\section{Main Results}\label{sec:2}
\subsection{Convergence rate of the density estimator}\label{lab:2}
In this section we obtain the asymptotic properties of the estimator \eqref{12}. To this end we derive the bias and the variance of the estimates in the following lemmas.
\begin{lemma}\label{lem1}If $b\rightarrow 0$ as $n\rightarrow \infty$, then the
bias of the pdf estimate  \eqref{12} is equal to
\begin{eqnarray}\label{21}Bias(\widehat{f}(x))&=& \left\{
\begin{array}{ll}
C_1(x,a)+hC_2(x,a,c_1)+o(h) &   \mbox{if}\qquad x\in
[0,a],
\\
B_1(x,a)+bB_2(x,a,c_2)+o(b)& \mbox{if}\qquad x> a.
\end{array}
\right.,
\end{eqnarray}
where we introduce the notations
\begin{eqnarray}\label{13}C_1(x,a)=\frac{xa}{2}f''(x),\quad C_2(x,a,c_1)\equiv c_1c_2(x,a)=c_1\left(f'(x)+f''(x)\frac{a}{2}+f'''(x)\frac{xa}{2}\right),
\end{eqnarray}
\begin{eqnarray}\label{14}&&B_1(x,a)=f\left(a\Gamma\left(t\right)\right)-f(x)+f''\left(a\Gamma\left(t\right)\right)\frac{a^2}{2}
\left(\Gamma\left(t\right)-\Gamma\left(r\right)\right)^2,
\end{eqnarray}
\begin{eqnarray}\label{15}&&B_2(x,a,c_2)\equiv c_2b_2(x,a)=\frac{a^2c_2}{x^2}\Bigg(-f'\left(a\Gamma\left(t\right)\right)\Gamma\left(t\right)\Psi\left(t\right)+f''\left(a\Gamma\left(t\right)\right)a\\\nonumber
&\cdot&\Bigg(\left(\Gamma\left(t\right)-\Gamma\left(r\right)\right)\left(\Gamma\left(t\right)\Psi\left(t\right)-2
\Gamma\left(r\right)\Psi\left(r\right)\right)\Bigg)-f'''\left(a\Gamma\left(t\right)\right)\frac{a^2}{2}\left(\Gamma\left(t\right)-\Gamma\left(r\right)\right)^2\Gamma\left(t\right)\Psi\left(t\right)\Bigg),
\end{eqnarray}
where we use the notations
\begin{eqnarray}\label{16}t&=&1+\frac{a}{x},\quad r = 1+\frac{2a}{x}.
\end{eqnarray}
\end{lemma}
\begin{lemma}\label{lem2}If $b\rightarrow 0$ as $n\rightarrow \infty$, then the
variance of the pdf estimate  \eqref{12} is equal to
\begin{eqnarray}\label{20}Var(\widehat{f}(x))&=&
\frac{1}{n}\Bigg(A_1(x,a)-(C_1(x,a)+f(x))^2+h(A_2(x,a,c_2)\\\nonumber
&-&2C_2(x,a,c_1)(C_1(x,a)+f(x)))\Bigg)+o(h)  \quad \mbox{if}\qquad x\in
[0,a],\nonumber\\\nonumber
Var(\widehat{f}(x))&=&
\frac{1}{n}\Bigg(D_1(x,a)-(B_1(x,a)-f(x))^2+b(D_2(x,a,c_2)\\\nonumber
&-&2B_2(x,a,c_2)(B_1(x,a)-f(x)))\Bigg)+o(b) \quad \mbox{if}\qquad x> a.\nonumber
\end{eqnarray}
where we introduce the notations
\begin{eqnarray}\label{18}&&A_1(x,a)=-f\left(x-\frac{a}{2}\right)\frac{\sqrt{x}}{\sqrt{a}(a-2x)},\\\nonumber
&&A_2(x,a,c_1)\equiv c_1a_2(x,a)=-c_1\Bigg(f\left(x-\frac{a}{2}\right)\frac{a+2x}{2\sqrt{ax}(a-2x)^2}\\\nonumber
&+&\frac{\sqrt{x}}{\sqrt{a}(a-2x)}\left(f'\left(x-\frac{a}{2}\right)+\frac{a}{4}\left(x-\frac{a}{2}\right)f''\left(x-\frac{a}{2}\right)\right)\Bigg)
\end{eqnarray}
and
\begin{eqnarray}\label{20}&&D_1(x,a)=\frac{x2^{\frac{3x}{a}-1}}{a^2}\left(f(2a)\left(\frac{x(x-3a)}{2a^2}+2\right)+f'(2a)(x-a)+f''(2a)2a^2\right),\\\nonumber
&&D_2(x,a,c_2)\equiv d_{21}(x,a)+c_2d_{22}(x,a)=\frac{x2^{\frac{3x}{a}-1}}{a^2}\Bigg(f(2a)\Bigg(\frac{(x-2a)(dx^2+a(c_2-dx))}{2a^3}\\\nonumber
&+&\frac{dx^2+a(c_2-dx)}{a}
-\frac{x(x-a)(x-2a)(6\gamma-10-\ln(4))}{2a^3}\Bigg)\\\nonumber
&+&f'(2a)\left(\frac{dx^2+a(c_2-dx)}{a}-2xd(x-2a+1)-\frac{x(x-a)(6\gamma-10-\ln(4))}{a}\right)\\\nonumber
&-&f''(2a)\left(2x(d(x-a)+2a^2+a(6\gamma-10-\ln(4)))\right)+\frac{c_22^{\frac{3x}{a}-1}}{a^3}\left(-x\ln(a)+a+2x\ln(2)\right)\\\nonumber
&\cdot&\left(f(2a)\left(\frac{x(x-3a)}{2a^2}+2\right)+f'(2a)(x-a)+f''(2a)2a^2\right)\Bigg).
\end{eqnarray}
\end{lemma}
The proofs of the latter lemmas are given in Appendices \ref{Ap1} and \ref{Ap2}.
\subsection{The optimal bandwidth parameters for the the density estimator}
To find the mean integrated squared error (MISE) we use the results of the two last paragraphs. Hence, we can write that for the domain  $x\in(0,a]$ the MSE is
\begin{eqnarray*}MSE(\hat{f}(x))_G&=&C_1^2(x,a)+h^2C_2^2(x,a,c_1)+2hC_1(x,a)C_2(x,a,c_1)\\\nonumber
&+&\frac{1}{n}\Bigg(A_1(x,a)-(C_1(x,a)+f(x))^2+h(A_2(x,a,c_1)\\\nonumber
&-&2C_2(x,a,c_1)(C_1(x,a)+f(x)))\Bigg)+o(h).
\end{eqnarray*}
Hence, from the minima of the latter equation we can write
\begin{eqnarray}\label{23}h_{opt}(x,a,n)&=&-\frac{C_1(x,a)}{C_2(x,a,c_1)}-\frac{1}{C_2(x,a,c_1)n}\left(\frac{A_2(x,a,c_1)}{2C_2(x,a,c_1)}-C_1(x,a)-f(x)\right).
\end{eqnarray}
Substituting the latter bandwidth in to the MISE we get the following rate
\begin{eqnarray*}MSE(\hat{f}(x))_{Gopt}&=&-\frac{1}{n^2}\left(\frac{A_2(x,a,c_1)}{2C_2(x,a,c_1)}-(C_1(x,a)+f(x))\right)^2\\
&+&\frac{1}{n}\left(A_1(x,a)-\frac{A_2(x,a,c_1)C_1(x,a)}{C_2(x,a,c_1)}+(C_1^2(x,a)-f^2(x))\right).
\end{eqnarray*}
For the  domain $x>a$ the MSE is
\begin{eqnarray*}MSE(\hat{f}(x))_W&=&B_1^2(x,a)+b^2B_2^2(x,a,c_2)+2bB_1(x,a)B_2(x,a,c_2)\\\nonumber
&+&\frac{1}{n}\Bigg(D_1(x,a)-(B_1(x,a)-f(x))^2+b(D_2(x,a,c_2)\\\nonumber
&-&2B_2(x,a,c_2)(B_1(x,a)-f(x)))\Bigg)+o(b)
\end{eqnarray*}
and the optimal bandwidth is
\begin{eqnarray}\label{24}b_{opt}(x,a,n)&=&\frac{-B_1(x,a)}{B_2(x,a,c_2)}-\frac{1}{B_2(x,a,c_2)n}\Bigg(\frac{D_2(x,a,c_2)}{2B_2(x,a,c_2)}-B_1(x,a)+f(x)\Bigg).
\end{eqnarray}
Substituting the latter bandwidth in to the MISE we get the following rate
\begin{eqnarray*}MSE(\hat{f}(x))_{Wopt}&=&-\frac{1}{n^2}\left(\frac{D_2(x,a,c_2)}{2B_2(x,a,c_2)}-(B_1(x,a)+f(x))\right)^2\\
&+&\frac{1}{n}\left(D_1(x,a)-\frac{D_2(x,a,c_2)B_1(x,a)}{B_2(x,a,c_1)}+(B_1^2(x,a)-f^2(x))\right).
\end{eqnarray*}
Since the condition $h_{opt}(a,a,n)=b_{opt}(a,a,n)c_2/c_1$ holds, we must find the parameters $a,c_1,c_2$ to satisfy the latter it.
Let us select the bandwidth $b_{opt}(a,a,n)$ which is optimal for the tail part of the estimate. Hence, the second bandwidth is $h_{b_{opt}}(a,a,n)=b_{opt}(a,a,n)c_2/c_1$. We can find such constants $a,c_1,c_2$  that
\begin{eqnarray*}&&\min\limits_{a,c_1,c_2}{h_{opt}(a,a,n)-h_{b_{opt}}(a,a,n)}
\end{eqnarray*}
holds.
Hence, substituting the values of the bandwidths we get the following condition
\begin{eqnarray}\label{22}&&c_2=\frac{1}{d_{21}(a,a)}\left(\frac{B_1(a,a)b_2(a,a)}{C_1(a,a)}
\left(\frac{a_2(a,a)}{c_2(a,a)}-2f(a)\right)-2f(a)b_2(a,a)-d_{22}(a,a)\right).
\end{eqnarray}
Hence, we can select any negative $c_1$, e.g. $c_1=-1$.
\begin{eqnarray*}d_{21}(x,a)&=&2^{\frac{3a}{2}}\Bigg(\frac{f'(2a)}{\ln(10)}(a-1)(\gamma\ln(10)-\ln(5))\\
&+&f''(2a)\left(\gamma-1+2a^2\left(\frac{\ln(4)}{\ln(10)}-6\gamma+10\right)+\frac{\ln(2)}{\ln(10)}\right)\Bigg),
\end{eqnarray*}
\begin{eqnarray*}d_{22}(x,a)&=&2^{\frac{3a}{2}-1}\Bigg(\frac{f(2a)}{a^2}\left(
\frac{2a-1}{2}+3\left(1-\frac{1}{\ln(10)(\ln(a)-2\ln(2))}\right)\right)\\
&+&\frac{f'(2a)}{a}+
f''(2a)\left(1-\frac{1}{\ln(10)}(\ln(a)-2\ln(2))\right)\Bigg).
\end{eqnarray*}
\section{Simulation study}
To investigate the performance of the Gamma-Weibull kernel estimator we select the following positive defined pdfs: the Weibull ($a =1, b=0.9$).
We generate Weibull  i.i.d samples with sample sizes $n\in\{100, 500, 1000 , 2000\}$ using standard Matlab generators.
\par Let us  find the value of the smoothing parameter using the rule
of thumb method. To this end we choose the gamma density
\begin{eqnarray} \label{9a}
f(x)=\frac{x^{\rho-1}\exp(-\frac{x}{\kappa})}{\kappa^{\rho}\Gamma(\rho)}
\end{eqnarray}
as a reference function.  Its first moment and the variance
are $\rho \kappa$ and $\rho \kappa^2$, respectively. According to the method
of moments, we have to equate them to the first sample moment $\bar
m = n^{-1}\sum_{i=1}^n X_i $ and the sample variance $\bar D =n^{-1}
\sum_{i=1}^{n}(X_i - \bar m)^2$, correspondingly. Then we obtain for
the parameters of \eqref{9a} following simple expressions
\begin{eqnarray}\label{10}
 \kappa_m = \bar D/\bar m,\quad \rho_m = (\bar m)^2/\bar D.
\end{eqnarray}
Hence, selecting some $a$ we can estimate the value of $c_2$ using \eqref{22}. Next,
calculating one of the optimal bandwidths \eqref{23} or \eqref{24} we can immediately find the other one using $h_{opt}(a,a,n)=b_{opt}(a,a,n)c_{2m}$.
\begin{eqnarray*}&&
\end{eqnarray*}

\section{Appendix}
\subsection{Proof of Lemma \ref{lem1}}\label{Ap1}
To find the bias of the estimate $\widehat{f}(x)$ let us write the expectation of the kernel estimator \eqref{12}
\begin{eqnarray}\label{17}E(\widehat{f}(x))&=& \left\{
\begin{array}{ll}
E_G(\widehat{f}(x))=\int_0^{\infty}K_{\rho(x,h),\theta}(y)f(y)dy = E(f(\xi_x)),&   \mbox{if}\qquad x\in
[0,a),
\\
E_W(\widehat{f}(x))=\int_0^{\infty}K_{k(x,b),\lambda}(y)f(y)dy = E(f(\eta_x)),& \mbox{if}\qquad x\geq a.
\end{array}
\right.
\end{eqnarray}
where $\xi_x$ is the gamma distributed $(\rho(x,h),\theta)$ r.v.s with the expectation $\mu_x=\rho(x,h)\theta$ and the variance $Var(\xi_x)=\rho(x,h)\theta^2$
and $\eta_x$ is the weibull distributed $(k(x,b),\lambda)$ r.v.s with the expectation $\widetilde{\mu}_x=\lambda\Gamma(1+\frac{1}{k(x,b)})$ and the variance $\widetilde{Var}(\eta_x)=\lambda^2\left(\Gamma(1+\frac{2}{k(x,b)})-\Gamma(1+\frac{1}{k(x,b)})\right)^2$.
\par Let us use the parameters \eqref{3} and $\theta=\lambda=a$. Hence, using the Taylor series in the point $\mu_x$ the expectation for the domain $x\in
[0,a]$ can be written as
\begin{eqnarray}\label{5}E(f(\xi_x))&=&f(\mu_x)+\frac{1}{2}f''(\mu_x)Var(\xi_x)+o(h)\nonumber\\
&=&f(x+c_1h)+\frac{a(x+c_1h)}{2}f''(x+c_1h)+o(h)\nonumber\\
&=&f(x)+f'(x)c_1h+\frac{a(x+c_1h)}{2}\left(f''(x)+f'''(x)c_1h\right)+o(h).
\end{eqnarray}
Thus, it is straightforward to verify that the bias of the estimate in the domain  $x\in
[0,a]$ is
\begin{eqnarray*}Bias_G(\widehat{f}(x))&=&C_1(x,a)+C_2(x,a)h+o(h),
\end{eqnarray*}
where we used the notations \eqref{13}.
\par To find the bias for the domain $x> a$ we need to Taylor expand $E(f(\eta_x))$ in the point $\widetilde{\mu}_x$. However the $\widetilde{\mu}_x$ and $\widetilde{Var}(\eta_x)$ contain the gamma function depending on the bandwidth parameter. To get there order by $b$ we need to expand them knowing that $b\rightarrow0$ and $nb\rightarrow\infty$ as the $n\rightarrow\infty$. Hence, we can write
\begin{eqnarray*}&&\widetilde{\mu}_x=a\Gamma\left(t\right)-b\frac{a^2c_2}{x^2}\Gamma\left(t\right)\Psi\left(t\right)+o(b),\\
&&\widetilde{Var}(\eta_x)=a^2\left(\Gamma\left(t\right)-\Gamma\left(r\right)\right)^2+b\frac{2a^3c_2}{x^2}\left(\Gamma\left(t\right)-\Gamma\left(r\right)\right)\\
&\cdot&\left(\Gamma\left(t\right)\Psi\left(t\right)-2\Gamma\left(r\right)\Psi\left(r\right)\right)+o(b),
\end{eqnarray*}
where we used the notations \eqref{16} and $\Psi(\cdot)$ is a dygamma function. Thus, the expectation \eqref{17} can be written as
\begin{eqnarray}\label{6}E(f(\eta_x))&=&f(\widetilde{\mu}_x)+\frac{1}{2}f''(\widetilde{\mu}_x)\widetilde{Var}(\eta_x)+o(h)\\
&=&f\left(a\Gamma\left(t\right)\right)-f'\left(a\Gamma\left(t\right)\right)\frac{a^2c_2\Gamma\left(t\right)\Psi\left(t\right)}{x^2}b\nonumber\\
&+&\frac{1}{2}\Bigg(a^2\left(\Gamma\left(t\right)-\Gamma\left(r\right)\right)^2
+\frac{2a^3c_2b}{x^2}\left(\Gamma\left(t\right)-\Gamma\left(r\right)\right)\left(\Gamma\left(t\right)\Psi\left(t\right)-2
\Gamma\left(r\right)\Psi\left(r\right)\right)\Bigg)\nonumber\\
&\cdot&\left(f''\left(a\Gamma\left(t\right)\right)
-f'''\left(a\Gamma\left(t\right)\right)\frac{a^2c_2b}{x^2}\Gamma\left(t\right)\Psi\left(t\right)\right)+o(b).\nonumber
\end{eqnarray}
Therefore, we can write that the bias of the pdf estimate in the domain $x>a$ is
\begin{eqnarray*}Bias_W(\widehat{f}(x))&=&B_1(x,a)+bB_2(x,a)+o(b),
\end{eqnarray*}
where we used the notations \eqref{14} and \eqref{15}.

\subsection{Proof of Lemma \ref{lem2}}\label{Ap2}
By definition the variance is
\begin{eqnarray}\label{4}
&& Var(\hat{f}(x))=\frac{1}{n}Var(K(x)) =\frac{1}{n}\left(E(K^2(x))-E^2(K(x))\right).
\end{eqnarray}
The second term of the right-hand side of \eqref{4} is the square of the
\eqref{5} and \eqref{6} for the domains $x\in[0,a]$ and $x>a$, respectively.
The first term of the right-hand side of \eqref{4} for the domain $x\in[0, a]$ can be
represented by
\begin{eqnarray}\label{8} E(K_G^2(x))=\int\limits_0^\infty K_G^2(y) f(y)dy=
\int\limits_0^\infty \frac{y^{2\left(\frac{x+c_1h}{a}-1\right)}e^{-2y/a}}{a^{2\frac{x+c_1h}{a}}\Gamma^2\left(\frac{x+c_1h}{a}\right)}f(y)dy=B(x,h,a)E(f(\zeta_x)),
\end{eqnarray}
where $\zeta_x$ is the gamma distributed with the parameters $\left(\frac{2(x+c_1h)}{a}-1,\frac{a}{2}\right)$ r.v.s with the expectation $\mu_{\zeta}=x+c_1h-\frac{a}{2}$ and the variance $Var(\zeta_x)=(x+c_1h)\frac{a}{2}-\frac{a^2}{4}$ and we used the following notation
\begin{eqnarray}\label{7}&&B(x,h,a)=\frac{\Gamma\left(\frac{2(x+c_1h)}{a}-1\right)}{a\Gamma^2\left(\frac{x+c_1h}{a}\right)2^{\frac{2(x+c_1h)}{a}-1}}.
\end{eqnarray}
Using the Stirling's formula for the gamma function and since $x\in(0,a]$ and $h\rightarrow 0$ as $n\rightarrow\infty$ we can expend \eqref{7} as
%\begin{eqnarray*}
%B(x,h,a)&=&\frac{2^{1-\frac{2x}{a}}\Gamma\left(\frac{2x}{a}-1\right)}{a\Gamma^2\left(\frac{x}{a}\right)}-
%h\frac{2^{2\left(1-\frac{x}{a}\right)}\Gamma\left(\frac{2x}{a}-1\right)}{a^2\Gamma^2\left(\frac{x}{a}\right)}\left(\Psi\left(\frac{a}{x}\right)-\Psi\left(t\right)+\ln 2\right)+o(h).
%\end{eqnarray*}
\begin{eqnarray*}
B(x,h,a)&=&-\frac{\sqrt{x}}{\sqrt{a}(a-2x)}-hc_1\frac{a+2x}{2\sqrt{ax}(a-2x)^2}+o(h).
\end{eqnarray*}
The expectation in \eqref{8} can be Taylor expanded similarly to the previous proof as
\begin{eqnarray*}E(f(\zeta_x))&=&f\left(x+c_1h-\frac{a}{2}\right)+\left((x+c_1h)\frac{a}{4}-\frac{a^2}{8}\right)f''\left(x+c_1h-\frac{a}{2}\right)+o(h)\\
&=&f\left(x-\frac{a}{2}\right)+c_1h\left(f'\left(x-\frac{a}{2}\right)+\frac{a}{4}\left(x-\frac{a}{2}\right)f''\left(x-\frac{a}{2}\right)\right)+o(h).
\end{eqnarray*}
Hence, the expectation \eqref{8} is
\begin{eqnarray*}&& E(K^2(x))=A_1(x,a)+hA_2(x,a)+o(h),
\end{eqnarray*}
where we used the notations \eqref{18}.
%\begin{eqnarray*}&&A_1(x,a)=f\left(x-\frac{a}{2}\right)\frac{2^{1-\frac{2x}{a}}\Gamma\left(\frac{2x}{a}-1\right)}{a\Gamma^2\left(\frac{x}{a}\right)},
%\end{eqnarray*}
%\begin{eqnarray*}&&A_2(x,a)=\frac{2^{2\left(1-\frac{x}{a}\right)}\Gamma\left(\frac{2x}{a}-1\right)}{a\Gamma^2\left(\frac{x}{a}\right)}\Bigg(
%\left(f'\left(x-\frac{a}{2}\right)+\frac{a}{4}\left(x-\frac{a}{2}\right)f''\left(x-\frac{a}{2}\right)\right)\\
%&-&f\left(x-\frac{a}{2}\right)\frac{2}{a}\left(\Psi\left(\frac{a}{x}\right)-\Psi\left(t\right)+\ln 2\right)\Bigg).
%\end{eqnarray*}
Hence, the variance \eqref{4} for the domain  $x\in(0,a]$ is
\begin{eqnarray*}&&Var_G(\hat{f}(x))=\frac{1}{n}\left(A_1(x,a)-C_1(x,a)+hc_1(A_2(x,a)-2C_1(x,a)C_2(x,a))\right)+o(h).
\end{eqnarray*}
For the domain $x>a$ we can write similarly to the previous part of the proof that
\begin{eqnarray}\label{9} E(K_W^2(x))&=&\int\limits_0^\infty K_W^2(y) f(y)dy=
\int\limits_0^\infty \frac{k(x,b)^2}{a^2}\left(\frac{y}{a}\right)^{2\left(k(x,b)-1\right)}\exp\left(-2\left(\frac{y}{a}\right)^{k(x,b)}\right)f(y)dy\nonumber\\
&=&\frac{4^{k(x,b)}k(x,b)}{a^{k(x,b)}}E(f(\varsigma_x)\varsigma_x^{k(x,b)-1}),
\end{eqnarray}
where $\varsigma_x$ is the weibull distributed r.v.s with the parameters $\left(k(x,b),2^{k(x,b)a}\right)$ and the expectation
\begin{eqnarray*}m_x=2(a-bxd)+o(b^2),\quad d= \gamma-1+\ln(2)\end{eqnarray*}
and the variance
\begin{eqnarray*}
Var_{m_x}&=&4a^2-4bax(6\gamma-10+\ln(4))+o(b^2),\end{eqnarray*}
where $\gamma$ is the Euler-Mascherson constant.
 Hence, the expectation \eqref{9} can be written as
\begin{eqnarray*}&&E(f(\varsigma_x)\varsigma_x^{k(x,b)-1})=f(m_x)m_x^{k(x,b)-1}+\frac{Var_{m_x}}{2} \Bigg(f''(m_x)m_x^{k(x,b)-1}\\
&+&2(k(x,b)-1)f'(m_x)m_x^{k(x,b)-2}+(k(x,b)-1)(k(x,b)-2)f(m_x)m_x^{k(x,b)-3}\Bigg)+o(b)\\
&=&m_x^{k(x,b)-1}\Bigg(f(m_x)+\frac{Var_{m_x}}{2} \Bigg(f''(m_x)+(k(x,b)-1)m_x^{-1}\Bigg(f'(m_x)+(k(x,b)-2)f(m_x)m_x^{-1}\Bigg)\Bigg)\Bigg).
\end{eqnarray*}
Using the Taylor series we can write that
\begin{eqnarray*}&&m_x^{k(x,b)-1}=(2a)^{\frac{x}{a}-1}\Bigg(1+\frac{b}{a}\Bigg(c_2\ln(2a)+xd(1-x)\Bigg)\Bigg)+o(b),\\
&&(k(x,b)-1)m_x^{-1}=\frac{x-a}{2a^2}+b\frac{a(c_2-dx)+dx^2}{2a^3}+o(b)\\
&&(k(x,b)-2)m_x^{-1}=\frac{x-2a}{2a^2}+b\frac{a(c_2-2dx)+dx^2}{2a^3}+o(b)\\
\end{eqnarray*}
\begin{eqnarray*}&&\frac{4^{k(x,b)}k(x,b)}{a^{k(x,b)}}=x4^{\frac{x}{a}}a^{-\frac{x}{a}-1}
+c_2b4^{\frac{x}{a}}a^{-\frac{x}{a}-2}(-x\ln(a)+a+x\ln(4))+o(b).
\end{eqnarray*}
Hence, the variance is the following
\begin{eqnarray*}&&Var_W(\hat{f}(x))=\frac{1}{n}\left(D_1(x,a)+bD_2(x,a,c_2)-\left(B_1(x,a)+bB_2(x,a)+f(x)\right)^2\right),
\end{eqnarray*}
where we used the notations
\begin{eqnarray*}&&
\end{eqnarray*}

\begin{thebibliography}{8}

% --- example of Journal article:
\bibitem{Aksoy}
 \textsc{Aksoy, H. } (2000). Use of Gamma Distribution in Hydrological Analysis.
 \textit{Turk J. Engin Environ Sci}, \textbf{24}, 419 -- 428.
\bibitem{Asmussen}
  \textsc{Asmussen, S. R.} (2003). Steady-State Properties of GI/G/1.  \textit{Applied Probability and Queues. Stochastic Modelling and Applied Probability}, \textbf{51}, 266--301.
 \bibitem{Dobrovidov:12}
 \textsc{Dobrovidov, A.V. \ \textup{and} Koshkin, G.M.  \ \textup{and}  Vasiliev, V. A.} (2012).
Non-parametric state space models.
\textit{Kendrick press}, USA.
 \bibitem{TaufikBouezmarnia:Rom}
\textsc{Bouezmarnia, T. \ \textup{and}  Rombouts, J.V.K.} (2007). Nonparametric density estimation for multivariate bounded data.
\textit{Journal of Statistical Planning and Inference}, \textbf{140}, 1, 139--152.

\bibitem{TaufikBouezmarnia:Rom2}
\textsc{Bouezmarnia, T. \ \textup{and}  Rombouts, J.V.K.} (2010). Nonparametric density estimation for positive times series.
\textit{Computational Statistics and  Data Analysis},
\textbf{54}, 2, 245--261.
\bibitem{BouSca}
  \textsc{Bouezmarnia, T. and Scaillet, O.} (2003).
Consistency of Asymmetric Kernel Density
Estimators and Smoothed Histograms
with Application to Income Data.
   \textit{Econometric
Theory}, \textbf{21}, 390--412.
\bibitem{Chen:20}
\textsc{Song Xi Chen} (2000).
     Probability density function estimation using gamma kernels.
    \textit{Annals of the Institute of Statistical Mathematics}
      \textbf{54}, 471--480.
\bibitem{Furman}
 \textsc{Furman, E.} (2008). On a multivariate Gamma distribution.
  \textit{Statist. Probab. Lett.}, \textbf{78}, 2353--2360.

\bibitem{Hurlimann}
 \textsc{H\"{u}rlimann, W.} (2001). Analytical Evaluation of Economic Risk Capital for Portfolios of Gamma Risks.  \textit{ASTIN Bulletin}, \textbf{31}, 107--122.
\bibitem{Mar_2015}\textsc{Markovich, L.A.} (2015). Nonparametric estimation of multivariate density and its derivative  by dependent data using gamma kernels.  (submitted in Journal of Nonparametric statistics) (arXiv:1410.2507)
\bibitem{Scailet}
  \textsc{Scaillet, O.} (2004).
Density Estimation Using Inverse and Reciprocal Inverse
Gaussian Kernels.
 \textit{Journal of Nonparametric Statistics},
 \textbf{16}, 217--226.

 \bibitem{Embrechts}
    \textsc{Embrechts, P. Kl\"{u}ppelberg, C. and Mikosch, T.} Modelling Extremal Events for Insurance and Finance Springer-Verlag, 648 pages,  \textbf{corr. 4th printing}, 1st ed. 1997, \textit{Springer}.
\end{thebibliography}
\end{document}